\documentclass[12pt]{article}
\usepackage{amssymb}
\usepackage{amsmath}
\usepackage{amsthm}
\usepackage{color}
\usepackage{comment}
\usepackage{geometry}		
\usepackage{graphicx}

\DeclareFontFamily{OT1}{pzc}{}
\DeclareFontShape{OT1}{pzc}{m}{it}{<-> s * [1.10] pzcmi7t}{}
\DeclareMathAlphabet{\mathpzc}{OT1}{pzc}{m}{it}

\let\originalleft\left
\let\originalright\right
\renewcommand{\left}{\mathopen{}\mathclose\bgroup\originalleft}
\renewcommand{\right}{\aftergroup\egroup\originalright}

\newtheorem{theorem}{Theorem}[section]

\newtheorem{corollary}[theorem]{Corollary}
\newtheorem{lemma}[theorem]{Lemma}

\theoremstyle{definition}
\newtheorem{definition}{Definition}[section]

\theoremstyle{remark}

\begin{document}

\title{Piecewise linear circle maps and\\ conjugation to rigid rational rotations
}
\author{Paul Glendinning, Siyuan Ma and James Montaldi  \\
\small{Department of Mathematics, University of Manchester,}\\
\small{Oxford Road, Manchester M13 9PL, UK}
}
\maketitle


\begin{abstract}
Criteria for piecewise linear circle homeomorphisms to be conjugate to a rigid rotation, $x\to x+\omega~({\rm mod}~1)$,  with rational rotation number $\omega$ are given. The consequences of the existence of such maps in families of maps is considered and the results are illustrated using two examples: Herman's classic family of piecewise linear maps with two linear components, and a map derived from geometric optics which has four components. These results show how results for piecewise smooth circle homeomorphisms with irrational rotation numbers have natural correspondences with the case of rational rotation numbers for piecewise linear maps. In natural families of maps the existence of a parameter value at which the map is conjugate to a rigid rotation implies linear scaling of the rotation number in a neighbourhood of the critical parameter value and no mode-locked intervals, in contrast to the behaviour of generic families of circle maps. 
\end{abstract}


\section{Introduction}\label{sect:intro}
There is a rich theory of maps of the circle which goes back to Poincar\'e and Denjoy in the late nineteenth and early twentieth centuries. For orientation preserving homeomorphisms the rotation number, an average rotation rate, is a topological invariant that determines many properties of the dynamics. If the rotation number is rational then there is at least one periodic orbit, whilst if the rotation number is irrational then there are no periodic orbits and, provided the map is sufficiently smooth ($C^2$ is enough) all orbits are dense in the circle \cite{Devaney1989}. Dynamics is often characterised by conjugacy classes. Two maps $f$ and $g$ are topologically conjugate if there is a homeomorphism $h$ from the domain of $f$ to the domain of $g$ such that $h\circ f=g\circ h$ (think of $h$ as a coordinate transformation). If $h$ is a diffeomorphism then the maps are differentiably (or $C^1$) conjugate. Results that go back to Denjoy state that a $C^2$ orientation preserving homeomorphism of the circle with irrational rotation number $\omega$ is topologically conjugate to the rigid rotation
\begin{equation}\label{eq:rigrot}
x_{n+1}=r_\omega (x_n)=x_n+\omega \ \ ({\rm mod}~1),
\end{equation}
see \cite{Khanin2003} for more historical background. Herman \cite{H1979} and others developed this theory and made more detailed statements about the smoothness classes of the conjugacy. For example \cite{Yoccoz1983} there is a set of irrational numbers $\mathcal{I}$ with full Lebesgue measure such that if $f$ is a $C^3$ orientation preserving diffeomorphism of the circle with rotation number $\omega\in \mathcal{I}$ then $f$ is differentiably conjugate to $r_\omega$. The set $\mathcal{I}$ is defined via Diophantine conditions on the rotation number. Herman \cite{H1979}  also considers some piecewise linear (PWL) models and establishes conditions under which they are conjugate to an irrational rigid rotation, and this theory has been further developed and generalized (e.g. \cite{Khanin2003,Liousse2005}). The aim of this paper is to see how these results carry over to the case of rational rotation number.

An orientation preserving homeomorphism $f$ of the circle, parametrized by $x\in [0,1)$, can be described via a lift $F:\mathbb{R}\to\mathbb{R}$, where $F$ is continuous, strictly monotonic, $F(x+1)=F(x)+1$ and $f(x)=F(x)~({\rm mod}~1)$. Two lifts of the same circle map differ by an integer \cite{Devaney1989}. Herman's example \cite{H1979}, see  section~\ref{sect:Hex}, has two points at which the gradient of the map changes. Such points are called break points, and for irrational rotation numbers these results have been generalised to piecewise smooth maps with an arbitrary number of break points \cite{AA2016, Adouani2012, Dzh2012, Liousse2005}. There is a natural jump condition on the derivatives of the map on either side of the break points which determines whether a piecewise smooth (PWS) map with a given irrational rotation number is conjugate to a rigid rotation, see section~\ref{sect:results} for more details. The aim of this paper is to describe how these results apply in the case of rational rotation numbers in the piecewise linear (PWL) case. Khmelev \cite{Khmelev2005} develops results for more general PWS systems, but these are of necessity weaker. In the piecewise linear case phenomena that are infinitely improbable in the more general context become possible, and in particular the existence of maps with dynamics conjugate to a rigid rational rotation can become much more prevalent (section~\ref{sect:family}, Lemma~\ref{lem:prev}). Indeed, in some simple examples, such as those in sections~\ref{sect:Hex} and \ref{sect:GMMex} it becomes a codimension one phenomenon.

When considering one parameter families of maps, rational rotation numbers are associated with mode-locking, i.e. intervals of the parameter on which the rotation number takes a constant rational value.  Mode-locking is generic at rational rotation numbers for natural families of circle homeomorphisms \cite{dMvS}.  

Our main results are given in slightly loose language below (precise statements are given in section~\ref{sect:results}, monotonic families are defined in Definition~\ref{def:monotonic}, the basic idea is that these are families on which the rotation number varies continuously and monoronically).
\begin{itemize}
\item A PWL orientation preserving circle homeomorphism is conjugate to a rigid rational rotation if and only if every break point is periodic. (Theorem~\ref{thm:equ}.)
\item If a PWL orientation preserving circle homeomorphism is conjugate to a rigid rational rotation then a natural jump condition on the slopes of the map holds. (Theorem~\ref{thm:Jcon}.)
\item If $\{f_\mu\}$ is a monotonic family of PWL orientation preserving circle homeomorphism and there is a parameter $\mu_c$ for which $f_{\mu_c}$ is conjugate to a rigid rational rotation, then there is no mode-locking region for that rotation number. (Theorem~\ref{thm:nomodelock}.)
\item If $\{f_\mu\}$ is a monotonic family of PWL orientation preserving circle homeomorphism and there is a parameter $\mu_c$ for which $f_{\mu_c}$ is conjugate to a rigid rational rotation, then the rotation number for parameters near $\mu_c$  scales linearly with the parameter.  (Theorem~\ref{thm:scaling}.)
\end{itemize}

We will illustrate these results with two examples: the classic model of Herman \cite{H1979} and a model derived from refraction in a periodic medium developed in \cite{GMM2024,Ma2024}.

The remainder of the paper is structured as follows. In section~\ref{sect:results} we introduce the notation and some of the classic results for PWL and PWS orientation preserving homeomorphisms. Our results are then stated more precisely than in the list above. In section~\ref{sect:conj} we prove the theorem in the first bullet point above and in section~\ref{sect:jump} we describe the jump condition. In sections~\ref{sect:acim} and ~\ref{sect:bounded} we prove analogues of results for PWS circle maps with irrational rotation number on the existence of absolutely continuous invariant measures and the growth of derivatives. In section~\ref{sect:family} we prove the first of our results for families of PWL circle homeomorphism and the scaling of rotation numbers as a function of the parameters is established in section~\ref{sect:scaling}. Sections~\ref{sect:Hex} and \ref{sect:GMMex} show how the results apply to two our examples and then section~\ref{sect:conc} contains a brief conclusion,
making the connection with pinching of mode-locked regions (also called shrinking points) in other piecewise smooth maps \cite{S2023}.

\section{Definitions and Statement of Results}\label{sect:results}
An orientation preserving homeomorphism $f$ of the circle with lift $F$ is piecewise smooth (PWS) if there exists a finite set (the break points) $\{b_i\}_1^n$, $0\le b_1<b_2<\dots <b_n<1$ such that $F$ is $C^2$ on each interval $(b_i,b_{i+1})$, $1\le i\le n-1$, and on $(b_n,b_1+1)$, and at each break point
\begin{equation}\label{eq:break}
F^\prime (b_{i,-})=\lim_{x\uparrow b_i}F^\prime (x)\ne \lim_{x\downarrow b_i}F^\prime (x)=F^\prime (b_{i,+}).  
\end{equation}
The lift is extended to $\mathbb{R}$ in the usual way by the relationship $F(x+1)=F(x)+1$. An orientation preserving homeomorphism $f$ is piecewise linear (PWL) if it is PWS and $F$ is affine on each interval $(b_i,b_{i+1})$, $1\le i\le n-1$, and on $(b_n,b_1+1)$. We will sometimes move between the map of the circle $f$ and its representation via a lift $F$ without comment. 

Note that the number of break points $n$ of a PWS orientation preserving homeomorphism of the circle has $n\ge 1$, whilst $n\ge 2$ for a PWL orientation preserving homeomorphism of the circle. Herman's classic example \cite{H1979}, to which we return in section~\ref{sect:Hex}, has two break points, and the example in section~\ref{sect:GMMex} has four. 

Every lift $F$ of an orientation preserving homeomorphism of the circle $f$ has a rotation number, defined by
\begin{equation}\label{eq:rotno}
\rho (F)=\lim_{m\to \infty}\frac{1}{m}\left( F^m(x)-x\right) .
\end{equation}
By standard theory \cite{Devaney1989} the limit defining $\rho (F)$ exists and is independent of $x$ and the rotation numbers of two lifts of the same map differ by an integer. It is therefore natural to define the rotation number of $f$ as $\rho(f) =\rho(F)~ {\rm mod}~1$. However, when considering families of maps we will usually work with a family of lifts so that the rotation number can vary continuously through integer values (there are other conventions that would also achieve this). For irrational 
$\rho (F)$ there are now a number of key results for PWS maps that echo Herman's original theory for smooth maps. These depend on jump conditions for the derivatives on each side of the break points. Let
\begin{equation}\label{eq:jumpb}
J(b_i)=\frac{F'(b_{i,+})}{F'(b_{i,-})}, \quad 1\le i\le n,
\end{equation}
(note that by definition $J(b_i)\ne 1$ and is independent of the choice of lift) and
\begin{equation}\label{eq:globalJ}
J(f)=\prod_1^n J(b_i).
\end{equation}
Recall that a measure $\nu$ on the circle is invariant for $f$ if $\nu (f^{-1}(A))=\nu (A)$ for all measurable subsets $A$ of the circle. It is an absolutely continuous invariant measure if, in addition, $\ell(A)=0$ implies that $\nu (A)=0$, where $\ell$ is Lebesgue measure on the circle (i.e. the Haar measure). The measure is a probability measure if the measure of the circle equals one. Clearly any measure for which the measure of the circle is finite can be rescaled to give a probability measure.  

Examples of the general results proved for PWS maps with breaks include
\begin{itemize}
\item if $\rho$ is irrational of bounded type (i.e. the coefficients of the continued fraction expansion of $\rho$ are bounded), $f$ is sufficiently smooth and $J(f)=1$ then the invariant probability measure of $f$ is absolutely continuous with respect to Lebesgue measure \cite[Theorem 1.2]{Dzh2012}. 
\item if $\rho$ is irrational, $f$ is sufficiently smooth and $J(f)\ne 1$ then the invariant probability measure of $f$ is singular with respect to the Lebesgue (more strictly, Haar) measure \cite[Theorem 1.4]{Dzh2012}. 
\end{itemize}
In the case of two break points more can be said \cite{Dzh2006,Dzh2018}. In the context of this paper it is the PWL results that are most relevant here. For this special case Liousse \cite{Liousse2005} has significantly stronger results. For PWL maps, $J(f)\equiv 1$ since $F'(b_{i,+})=F'(b_{i+1,-})>0$ taking the indices so that $b_{n+1}=b_1+1$. This implies that these terms cancel in the product (\ref{eq:globalJ}).  The possibility of similar cancellations on a smaller set of indices is used in \cite{Liousse2005}: if there exists $p$, $1< p<n$, and and a set $S$ containing $p$ distinct indices in $\{1, \dots ,n\}$ such that
\begin{equation}\label{eq:trivialJ}
\prod_{i\in S} J(b_{i})=1,
\end{equation}
then we say that $f$ has trivial cancellations on $S$ (Liousse uses the word `\emph{compensations}' in French \cite{Liousse2005}).

An irrational number is of bounded type if the coefficients of its continued fraction expansion are bounded.

\begin{theorem}\label{thm:liousse2005}(\cite{Adouani2012,Liousse2005}) Suppose that $f$ is a PWL orientation preserving homeomorphism of the circle with $n$ break points, $n\ge 2$ and has rotation number $\rho$ which is irrational of bounded type. The following four statements are equivalent.
\begin{itemize}
\item[(i)] The union of the forwards and backwards orbit of every break point contains another break point and if $S$ contains the set of indices of break points on the same orbit then $f$ has trivial cancellations on $S$.
\item[(ii)] $f$ has an absolutely continuous invariant probability measure.
\item[(iii)] $f$ is piecewise $C^1$ conjugate to the rigid rotation $r_\rho$.
\item[(iv)] The number of break points of $f^n$ is bounded by $C$ for all $n$, with $C$ independent of $n$.
\end{itemize}
\end{theorem}

The equivalence of (i) and (ii) combines Proposition 1, Theorem 1 and Theorem 2 case (2) of \cite{Liousse2005}.  The equivalence of (ii) and (iv) is case 1 of Theorem 2 of \cite{Liousse2005}. The equivalence of cases (i) and (iii) follows from \cite[Corollary 1.6]{Adouani2012} (note that this part of the equivalence does not require $\rho$ to be of bounded type).

The simplifications associated with PWL maps makes it possible to prove rather more for PWL maps with rational rotation number than is possible for the PWS maps considered in \cite{Khmelev2005}. Our main results can now be written down precisely.

The first result shows that the periodicity of the break points, which is an obvious consequence of conjugacy, is actually also sufficient to imply conjugacy. This is a natural extension of the equivalence of Theorem~\ref{thm:liousse2005}(i)  to the case of rational rotation numbers.

\begin{theorem}\label{thm:equ}Let $f$ be a PWL orientation preserving homeomorphism of the circle with rational rotation number $p/q$. The map $f$ is conjugate to the rigid rational rotation $r_{p/q}$ of (\ref{eq:rigrot}) by a PWL homeomorphism if and only if every break point is periodic.\end{theorem}

Of course, if an orientation preserving homeomorphism of the circle has distinct periodic orbits then the orbits must have the same period, so the fact that the break points lie on periodic orbits with the same period is implicit in the statement. 

The next result describes the nature of the break point orbits and completes the equivalence with Theorem~\ref{thm:liousse2005}(i), but note that this is a consequence of Theorem~\ref{thm:equ} not an additional criterion.

\begin{theorem}\label{thm:Jcon}Let $f$ be a PWL orientation preserving homeomorphism of the circle with $n$ break points and rational rotation number. If $f$ is topologically conjugate to a rigid rational rotation then $\{1, 2, \dots ,n\}$ can be partitioned into $K$, $1\le K \le n/2$, disjoint sets $S_r$, $r=1, \dots ,K$ such that each set $S_r$ contains at least two elements and if $u,v\in S_r$ then $b_u$ and $b_v$ lie on the same periodic orbit and such that $f$ has trivial cancellations on $S_r$. 
\end{theorem}

Some simple corollaries are mentioned in section~\ref{sect:jump} and section~\ref{sect:acim}. Bringing these results together gives the first three parts of the equivalent result to Theorem~\ref{thm:liousse2005}. The fourth and fifth parts are proved in section~\ref{sect:bounded}.

\begin{theorem}\label{thm:ourequiv}Suppose that $f$ is a PWL homeomorphism of the circle with $n$ break points, $n\ge 2$, and has rational rotation number $\rho=p/q$. The following five statements are equivalent.
\begin{itemize}
\item[(i)] The union of the forwards and backwards orbit of every break point contains at least one other break point and if $S$ contains the set of indices of break points on the same orbit then $f$ has trivial cancellations on $S$.
\item[(ii)] $f$ has an absolutely continuous invariant probability measure with support equal to $\mathbb{S}^1$.
\item[(iii)] $f$ is conjugate to a rigid rotation by a PWL conjugating function.
\item[(iv)] The number of break points of $f^k$ is bounded by $N_1$ for all $k\ge 1$, with $N_1$ independent of $k$.
\item[(v)] $(f^k)^\prime (x_\pm)$ is bounded by $N_2$ for all $k\ge 1$, with $N_2$ independent of $k$.
\end{itemize}
\end{theorem}

In the next two results we will consider families of orientation preserving homeomorphisms, $\{f_\mu\}$, for $\mu$ in some open interval. 

\begin{definition}Consider a family $\{f_\mu\}$ of orientation preserving circle homeomorphisms with $\mu\in (a,b)$. $\{f_\mu\}$ is a continuous family if there is a family of lifts $\{F_\mu\}$ such that $F_\mu (x)$ is a continuous function of $\mu$ at fixed $x$ and continuous in $x$ at fixed $\mu$. The family is monotonic if it is continuous and  $\mu_1<\mu_2$ implies that $F_{\mu_1}(x)<F_{\mu_2}(x)$ for all $x\in\mathbb{R}$. 
\label{def:monotonic}\end{definition}

Note that the definition (\ref{eq:rotno}) implies that for a monotonic family the rotation number is a continuous function of $\mu$, e.g. \cite{RT1991}, and for a monotonic family then $\rho (F_{\mu_1})\le\rho (F_{\mu_2})$. We can deal with decreasing rotation numbers by reversing the sign of $\mu$.  A nice way to construct monotonic families of circle homeomorphisms is to take the lift $G(x)$ of an orientation preserving circle homeomorphism and consider the family of circle maps with lifts $\mu  +G(x)$.  

For a monotonic family of maps $\{f_\mu\}$, $a<\mu <b$, let
\[
M_\alpha = \{\mu ~|~ \rho (F_\mu)=\alpha \},
\]
and suppose that $\alpha\in (\rho (F_a),\rho (F_b))$.   
so if $M_\alpha$ is a non-trivial interval then mode-locking occurs. For families of PWL circle homeomorphisms is that there is no mode-locking at parameter values for which the map is conjugate to a rigid rational rotation.

\begin{theorem}\label{thm:nomodelock}Let $\{f_\mu\}$ be a monotonic family of PWL orientation preserving homeomorphisms of the circle. If there exists $\mu_c$ such that $f_{\mu_c}$ is conjugate to a rigid rational rotation with rotation number $p/q$ then $M_{p/q}=\{\mu_c\}.$ 
\end{theorem}

This is not a new result, though we include the proof here for convenience. It is essentially the content of \cite[Ch.~1, Lemma 4.1]{dMvS} for example, where it is shown that if $f_{\mu}^q$ is not the identity then $M_{p/q}$ is a non-trivial interval. 

\begin{figure}[t!]
\centering
\includegraphics[height=5.5cm]{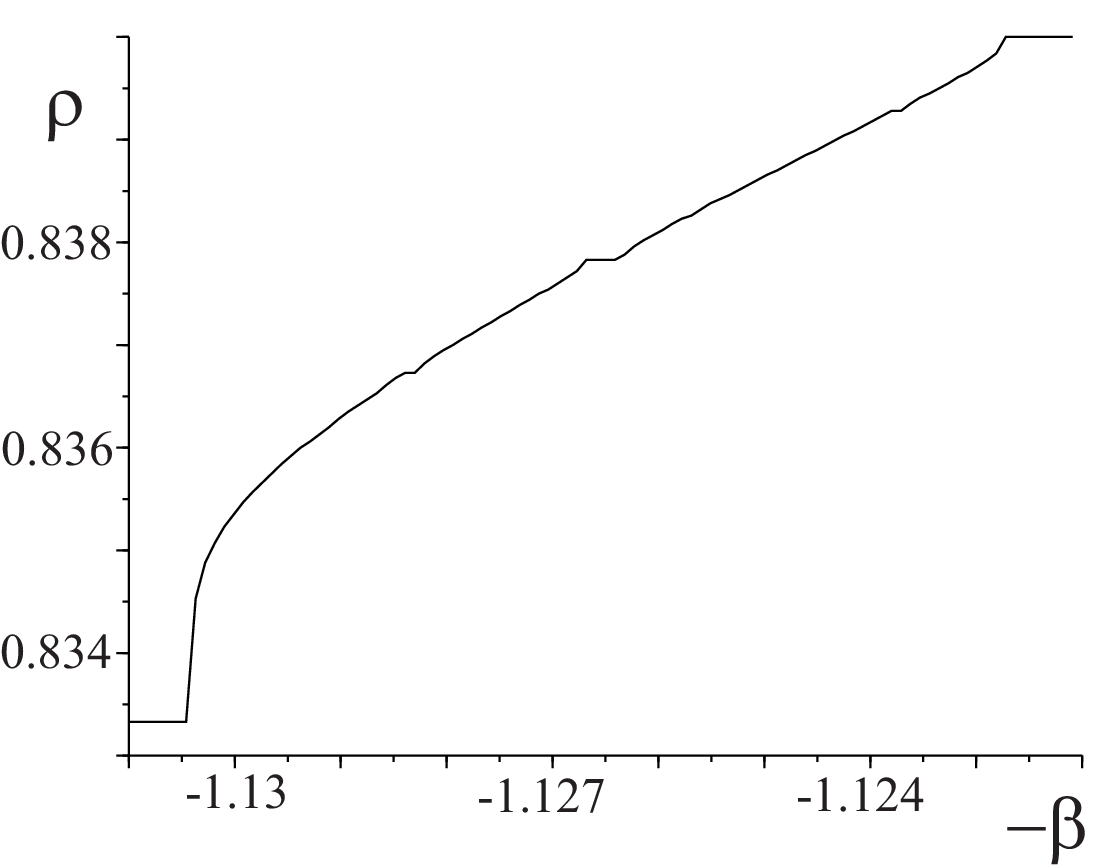}
\caption{
Plot of the rotation number $\rho$ (vertical axis) of (\ref{eq:GMMex}) against the parameter $-\beta$. The rotation number is calculated using 100000 iterates of the map at 100 equidistant values of the parameter with $\alpha=2$ and $1.122\le \beta < 1.131$. Part of the mode locked region with rotation number $\frac{5}{6}$ is visible on the left of the figure. The parameter $-\beta$ is used in the plot so that the rotation number increases with parameter as is the convention in the theoretical part of this paper.   
}
\label{fig:lockingintro}
\end{figure}

Our final result concerns the scaling of the rotation number. In the generic case there is no conjugacy to a rigid rational rotation (see section~\ref{sect:family}) and for one-parameter families rational rotation numbers occur on intervals of parameters. These intervals provide the stairs of the `Devil's staircase' when the rotation number is plotted against the parameter.  At the boundaries of mode locked regions the rotation number increases (or decreases) rapidly. This typical behaviour is shown in Figure~\ref{fig:lockingintro}. This figure shows the rotation number as a function of a parameter for the map of section~\ref{sect:GMMex} in a neighbourhood one end point of the mode locked region with rotation number $\frac{5}{6}$. The rotation number is constant on an interval of the parameter, the mode locked region, and then changes rapidly at each end, though there are evidently further plateaus of mode locked regions on the curve. In particular the modulus of slopes of bounding envelopes of the rotation number tends to infinity as the ends of the mode locked regions are approached from outside the model locked region. This typical behaviour has been described in detail for circle  maps which have intervals on which the map is constant \cite{Carretero1997}. The important difference for PWL families which have parameters where the map is topologically conjugate to a rigid rotation is that unlike the generic case of Figure~\ref{fig:lockingintro} the local scaling of the rotation number is linear, see Theorem~\ref{thm:scaling} below.

Families of PWL circle maps with $n$ break points are naturally described by $n$ pairs of continuous functions $(b_k(\mu ), \phi_k(\mu))_1^n$ where $b_k$ denotes the break points and $\phi_k$ the value of the function at the $k^{th}$ break point, i.e. if $F_\mu$ is the lift of the map then there is a neighbourhood in parameter space such that the break points are 
\begin{equation}\label{eq:bk}
b_1(\mu)<b_2(\mu )< \dots <b_{n}(\mu)<b_{n+1}=b_1(\mu )+1,
\end{equation}
and if 
\begin{equation}\label{eq:phik}
F_\mu (b_k(\mu ))=\phi_k(\mu ), 
\end{equation}
then the slope of the map, $s_k$ between $b_k$ and $b_{k+1}$ is
\begin{equation}\label{eq:sk}
s_k(\mu )=\frac{\phi_{k+1}-\phi_{k}}{b_{k+1}-b_{k}}. \quad s_k(\mu)\ne s_{k-1}(\mu ) .
\end{equation}
Thus for $k=1, \dots ,n$,
\begin{equation}\label{eq:Fk}
F_\mu (x)=\phi_k(\mu) +s_k(\mu )(x-b_k(\mu ))\ \ \ {\rm if}\ x\in [b_k(\mu),b_{k+1}(\mu)).
\end{equation}

If $\{f_\mu\}$ is a family of orientation preserving homeomorphisms and $f_0$ is conjugate to a rigid
rational rotation then the orbit under $f_0$ of each break point contains at least one other break point (Theorem~\ref{thm:Jcon}). Considering perturbations of such maps the images and preimages of these maps need to be controlled. This occupies most of the proof of Theorem~\ref{thm:scaling} in sections ~\ref{sect:family} and \ref{sect:scaling}. In the statement below the transversality conditions are given implicitly. 
      
\begin{theorem}\label{thm:scaling}Let $\{f_\mu\}$ be a family of PWL orientation preserving homeomorphisms, defined by $C^2$ functions $b_k(\mu )$ and $\phi_k(\mu )$ of (\ref{eq:bk}) and (\ref{eq:phik}), $1\le k\le n$. If there exists $\mu_c$ such that $f_{\mu_c}$ is conjugate to a rigid rational rotation with rotation number $p/q$ and 
the transversality condition (\ref{eqlem:linint}) holds then there exists an neighbourhood ${\mathcal N}$ of $\mu =\mu_c$ and constants $R_1>0$ and $R_2>0$ such that
\begin{equation}\label{eq:sc2}
|\rho (F_\mu )-\frac{p}{q}-R_1(\mu -\mu_c)|<R_2(\mu-\mu_c)^2.
\end{equation}
for all $\mu$ in ${\mathcal N}$. 
\end{theorem}

The transversality condition (\ref{eqlem:linint}) ensures that the natural 
parametrization for changes in the maps is order $(\mu -\mu_c)$ rather than (for example) a power of this difference.  Of course, 
another way to see (\ref{eq:sc2}) is that $\rho(\mu )$ is differentiable at $\mu_c$  with derivative $R_1$. The value of $R_1$ can be given explicitly in terms of properties of $F_{\mu_c}^q$:
\begin{equation}
R_1=\frac{1}{q\sum \kappa_i}
\label{eq:R1}
\end{equation}
where the sum is over the piecewise linear components of $F^q_{\mu_c}$ in $[0,1)$ and the coefficients $\kappa_i$ are given in the statement of Lemma~\ref{lem:timei}. 

\section{Conjugacy and dynamics of break points}\label{sect:conj}
This section contains the proof of Theorem~\ref{thm:equ}. This will need a preparatory lemma.

\begin{lemma}\label{lem:conjid}Suppose that $f$ is a PWL orientation preserving homeomorphism of the circle with rotation number $p/q$. If $f^q\equiv Id$ then $f$ is conjugate to a rigid rotation with rotation number $p/q$. The conjugacy can be chosen to be piecewise $C^k$ for any $k\ge 1$, and even PWL.
\end{lemma}

\emph{Proof:} Let $b_1$ be the first break point (indeed, any break point) and let $f^k(b_1)$ be the next element on the orbit of $b_1$ measured anti-clockwise. Note that the well-ordered property of orbits implies that if $K\ge 2$ in Theorem~\ref{thm:Jcon} then each of the other orbits of break points has one and only one point in the arc $(b_1, f^k(b_1))$ \cite{MT1984}. 

Let $h:[b_1,f^k(b_1)]\to [0, 1/q]$ be an arbitrary $C^k$ diffeomorphism (resp. an affine function) with $h(0)=0$ and $h(f^k(b_1))=1/q$. Now extend $h$ to $[f(b_1),f^{k+1}(b_1)]$ by $h(y)=h(f^{-1}(y))+p/q$, $y\in [f(b_1),f^{k+1}(b_1)]$.
Note that since $f$ is affine on $[b_1,f^k(b_1)]$, $h$ is $C^k$ (resp. affine) on  $(f(b_1),f^{k+1}(b_1))$.  Since $y=f(x)$ for some $x\in [b_1,f^k(b_1)]$ this implies that $h$ is well-defined on $[f(b_1),f^{k+1}(b_1)]$ and if $x\in [b+1,f^k(b_1)]$ then
\[
h\circ f(x)=r_{p/q}\circ h(x),
\]
where $r_{p/q}$ is the rigid rotaion (\ref{eq:rigrot}). Continuing this process on the iterates of $[b_1,f^k(b_1)]$, $h$ extends to the entire circle and the consistency condition on the final iteration (i.e. that the induced form of $h$ on $[b_1,f^{k}(b_1]$ from $[f^{-1}(b_1),f^{k-1}(b_1)]$ is the same as the $h$ defined at the beginning of the process follows directly from the assumption that $f^q\equiv Id$. Continuity at the break points follows from the continuity of the maps and periodicity of the break points.
\newline\rightline{$\square$}
 
\medskip
\emph{Proof of Theorem~\ref{thm:equ}}

Conjugacy to a rigid rotation with rotation number $p/q$ clearly implies that all points are periodic of period $q$, and in particular the break points are periodic with period $q$.

Now suppose that every break point is periodic. Since $f$ is a homeomorphism they must all have the same period, $q$, and denote the number of distinct orbits by $K$ . Let $p_i$, $1\le i\le qK$ denote the distinct points on the orbits of the break points each of period $q$, ordered so that
$0\le p_1<p_2< \dots <p_{qK}<1$. Note that $f^r$ is affine on $(p_i,p_{i+1})$ (with indices such that $p_{qK+1}=p_1$) for all $r$ as $(p_i,p_{i+1})$ contains no pre-images of break points. 

By definition $f^q(p_i)=p_i$, and so since $f^q$ is affine on each of the intervals it is the identity on each of the intervals, and hence on the circle. The result now follows from Lemma~\ref{lem:conjid}.
\newline\rightline{$\square$}

\section{Conjugacy and jump conditions}\label{sect:jump}

\emph{Proof of Theorem~\ref{thm:Jcon}:}
Since $f$ is conjugate to a rigid rational rotation with rotation number $\frac{p}{q}$, every point is periodic of period $q$ and there is a lift, $F$  of $f$ such that $F^q(x)=x+p$. The derivatives of $F^q$ from above and below exist and are equal if a point is not on the orbit of a break point. 

Suppose that $b_r$ is a break point and there is no other break point on the orbit. Then the derivatives from below and above of $F^q$ are
\[
F'(b_{k,-})\prod_1^{q-1}F'(F^r(b_k)) \quad {\rm and} \quad F'(b_{k,+})\prod_1^{q-1}F'(F^r(b_k))
\]
respectively. Since $b_k$ is the only break point on the orbit, the products do not need upper and lower limits, and since the slope of $F^q$ is one everywhere, this implies the two products of derivatives are equal (and equal to one) and hence that $F'(b_{k,-})=F'(b_{k,+})$, contradicting the assumption that $b_k$ is a break point.

Hence the orbit of every break point contains at least one other break point, and the total number of disjoint orbits of the $n$ break points is less than or equal to $n/2$.

By a similar argument, if an orbit of a break point contains precisely $r$ break points $b_{j_1}, \dots , b_{j_r}$ then the derivatives of $F^q$ from below and above at any of these break points are  
\[
C\prod_1^{r}F'(b_{j_r,-}) \quad {\rm and} \quad C\prod_1^{r}F'(b_{j_r,+}),
\]      
where $C$ is the product of the derivatives on points of the orbit that are not break points. Hence
\begin{equation}\label{eq:eqprod}
\prod_1^{r}F'(b_{j_r,-}) =\prod_1^{r}F'(b_{j_r,+}),
\end{equation}
and dividing through by the left hand side gives the jump condition (\ref{eq:trivialJ}) implying that $f$ has trivial cancellations on this orbit as stated in Theorem~\ref{thm:Jcon}.
\newline\rightline{$\square$}

\section{Absolutely continuous invariant measures}\label{sect:acim}
Theorem~\ref{thm:equ} has a simple corollary in terms of the existence of an absolutely continuous probability invariant measure (acipm).

\begin{corollary}\label{cor:haar}If a PWL homeomorphism of the circle $f$ is conjugate to a rigid rational rotation with rational rotation number $\frac{p}{q}$ then $f$ has an absolutely continuous invariant probability measure with support equal to  $\mathbb{S}^1$ .\end{corollary}

\emph{Proof:} Clearly $f^q$ preserves Lebesgue (Haar) measure as it is the identity. The result now follows from \cite{dMvS} (Chapter V, Lemma 3.1: if an induced map has an absolutely continuous invariant measure then so does the original map). More explicitly, if $\nu$ is an absolutely continuous invariant probability measure (acipm), define $f_*\nu$ by $f_*\nu (A)=\nu(f^{-1}(A))$ for any measurable set $A$. Since $f^q$ preserves Lebesgue measure $\ell$, if 
\[ \mu=\sum_0^{q-1}f_*^k\ell,
\]
then
\[
f_*\mu = \sum_1^qf_*^k\ell= f_*^q\ell+\sum_1^{q-1} f_*^k\ell,
\]
and since $f^q_*\ell=\ell $ as $f^q=Id$, $f_*\mu =\mu$ so $\mu$ is an invariant measure for $f$. Moreover, clearly if $\ell(A)=0$ then $\mu (A)=0$ so $\mu$ is absolutely continuous with respect to $\ell$. 

Since the support of $\ell$ is the whole circle then clearly the support of $\nu$ is also  $\mathbb{S}^1$.
\newline\rightline{$\square$}

The proof of the equivalence of (i) and (ii) in Theorem~\ref{thm:ourequiv} is completed by the following lemma. 

\begin{lemma}\label{lem:sing}If a PWL homeomorphism of the circle $f$ with rational rotation number $p/q$ is not conjugate to a rigid rational rotation then $f$ has no absolutely continuous invariant probability measure $\nu$ with ${\rm supp}\,(\nu )=\mathbb{S}$.\end{lemma}

\emph{Proof:}   If $f^q$ is not the identity then there exists $x\in \mathbb{S}^1$ such that $f^q(x)\ne x$. So there is a maximal open neighbourhood $U=(u_1,u_2)$ of $x$ such that $f^q(y)\ne y$ for all $y$ in $U$ and $f^q(u_i)=u_i$, $i=1,2$. For any $y\in U$, either $f^{-nq}(y) \to u_1$ or $f^{-nq}(y)\to u_2$. Thus $J_n=(f-{(n+1)q}(y),f^{-nq}(y))$ are disjoint arcs in $U$ (defined clockwise or counterclockwise depending on the limit) and if $nu$ is an acipm with ${\rm supp}\,(\nu )=\mathbb{S}$,  then $\nu (J_n)>0$ and $\nu (J){n+1}=\nu (J_n)$, a contradiction as $\sum \nu (J_n)<\infty$. 
\newline\rightline{$\square$}

Note that it is not hard to construct maps with rational rotation number and which are not conjugate to a rigid rotation yet have an acipm whose support is not the whole circle. For example consider the map which is the identity on a disjoint finite union of closed intervals and a piecewise linear interpolant between successive end-points on which $f$ is greater than $x$. Then the rotation number is zero, the map has an acipm (Lebesgue measure restricted to the intervals) and is not conjugate to a rigid rotation (not all points are periodic).

\section{Growth results}\label{sect:bounded}
In this section we establish the equivalence of the conjugacy to a rigid rotation to the bounded growth conditions of parts (iv) and (v) of Theorem~\ref{thm:ourequiv}. Note that we may assume Theorem~\ref{thm:equ}, proved in section~\ref{sect:conj}.

\begin{lemma}\label{lem:bddbreakpts}Suppose that $f$ is a $PWL$ circle homeomorphism with rotation number $p/q$. Then $f$ is conjugate to a rigid rotation if and only if the number of break points of $f^k$ is bounded by $qK$ where $K$ is the number of distinct orbits of the break points of $f$, for all $k\ge 1$. In particular $K$ is less than the number of break points of $f$.\end{lemma} 

\emph{Proof:} If $f$ is conjugate to a rigid rational rotation, each break point of $f$, $\{b_j\}_1^n$ lies on a periodic orbit, necessarily of period $q$, (Theorem~\ref{thm:equ}) and there are $K\le n/2$ of these (Theorem~\ref{thm:Jcon}). Thus the set of preimages of the breakpoints contains at most $qK$ points, and since any breakpoint of $f^k$ must be a preimage of a breakpoint of $f$ the number of break points of $f^k$ is bounded by $qK$.

If the number of break points of $f^k$ is bounded then each break point must lie on a periodic orbit (otherwise it would tend to a periodic orbit in backwards time and hence $\{f^{-r}(b_j)\}_{r\ge 0}$ would be infinite for some $j$ and since the rotation number is $p/q$ each of these is period $q$). Hence the map is conjugate to a rigid rotation by Theorem~\ref{thm:equ}.
\newline\rightline{$\square$}

\begin{lemma}\label{lem:bddslope}Suppose that $f$ is a $PWL$ circle homeomorphism with rotation number $p/q$. Then $f$ is conjugate to a rigid rotation if and only if the number of break points of $f^k$ $(k\ge 1$) is bounded by $s_M^{q-1}$ where $s_M$ is the the largest slope of the piecewise linear components of $f$.\end{lemma}

\emph{Proof:} If $f$, with lift $F$,  is conjugate to a rigid rotation then $(F^q)^\prime (x)=1$ for all $x$. Hence if $k=mq+r$, $0\le r<q$,
\[      
(F^k)^\prime (x_\pm )=(F^{mq})^\prime (F^r(x)_\pm)(F^r)^\prime (x_\pm)=(F^r)^\prime (x_\pm)\le s_M^{q-1}.
\]

Conversely, suppose that $f$ is not conjugate to a rigid rational rotation. Then there is a lift $F$ such that there is $z\in\mathbb{R}$ with $F^q(z)=z+p$ and $y\in \mathbb{R}$ such that $F^q(y)\ne y+p$. Suppose $F^q(y)<y$. Let $z$ be the smallest periodic point greater than $y$. On the PWL component of $F^q$ to the left of $z$, $(F^q)^\prime>1$ and hence $(F^{qm})^\prime (z_-)=[(F^q)^\prime (z_-)]^m\to \infty$ as $m\to \infty$ so the derivatives are not bounded. An analogous argument holds if $F^q(y)>y+p$.
\newline
\rightline{$\square$}

\section{Families of maps and the no mode-locking theorem}\label{sect:family}
As noted in section~\ref{sect:results} the rotation number of a continuous family of PWL circle maps is a continuous function of the parameter and if the family is monotonic then the rotation number is a monotonic function of the parameter.

It is natural to ask about the prevalence of parameters at which a map is conjugate to a rigid rational rotation in families of PWL maps. Lemma~\ref{lem:prev} shows that this is a codimension $n$ phenomenon, where $n$ is the number of break points.

\begin{lemma}\label{lem:prev}Suppose that $\{f_\mu\}$ is a generic family of PWL orientation preserving circle homeomorphisms with $n$ break points. The existence of a parameter at which the map is conjugate to a rigid rotation requires $n$ relationships to be satisfied simultaneously. \end{lemma}

\emph{Proof:} Since $f_\mu$ is PWL, Theorem~\ref{thm:Jcon} implies that the break points lie on the union of $K$, $1\le K\le n/2$ periodic orbits each of which contains at least two break points. Given $K$ points on these orbits this implies there are $n-K$ condition to ensure that the remaining break points are on the corresponding orbits, and one extra condition (i.e. $K$ in total)  on each to ensure that these orbits are actually periodic.
\newline\rightline{$\square$}   

\begin{corollary}If $\{f_\mu\}$ is a generic family of PWL orientation preserving circle homeomorphisms with two break points, then the existence of parameters at which the dynamics is conjugate to a rigid rational rotation is a codimension two phenomenon.\end{corollary}
 
Although this suggests that the existence of parameter values at which a map is conjugate to a rigid rational rotation are generally rare in maps with multiple break points (at least codimension $n$), the existence of symmetries can reduce the number of independent conditions in Lemma~\ref{lem:prev}, as is the case for the example from \cite{GMM2024,Ma2024} which will be discussed in section~\ref{sect:GMMex}.  

Theorem~\ref{thm:nomodelock} shows that for monotonic families of PWL circle homeomorphisms with break points there is no mode-locking at parameters at which the map is conjugate to a rigid rotation.

\medskip
\emph{Proof of Theorem~\ref{thm:nomodelock}:} Suppose that there is a parameter $\mu_c$ at which the map is conjugate to a rigid rotation with rotation number $p/q$. Then there is a lift of $f$ such that $F^q_{\mu_c}(x)=x+p$ for all $x\in\mathbb{R}$. Since $F_\mu$ is a monotonic family this implies that there is an open neighbourhood $U$ of $\mu_c$ such that $F^q_\mu (x)<x+p$ for all $x\in\mathbb{R}$ if $\mu<\mu_c$, $\mu\in U$, and  $F^q_\mu (x)>x+p$ for all $x\in\mathbb{R}$ if $\mu>\mu_c$, $\mu\in U$. In particular $\rho(F_\mu)\ne p/q$ if $\mu\in U\backslash\{\mu_c\}$ (there are no orbits of period $q$) and if $\mu \in U$ then $\rho (F_\mu)<p/q$ if $\mu<\mu_c$ and $\rho (F_\mu)>p/q$ if $\mu>\mu_c$.
\newline\rightline{$\square$}

\medskip
Recall from section~\ref{sect:results} that families of PWL circle maps with $n$ break points are naturally described by $n$ pairs of continuous functions $(b_k(\mu ), \phi_k(\mu))_1^n$ where $b_k$ denotes the break points and $\phi_k$ the value of the function at the $k^{th}$ break point, see (\ref{eq:bk}) and (\ref{eq:phik}). The slope of the map between $b_k$ and $b_{k+1}$ is $s_k$ defined by (\ref{eq:sk}). This formulation can also be used to identify monotonic families. A brute force calculation of
the difference between the values of the map at different parameter values (see the proof of Lemma~\ref{lem:linint} below) shows that  the lifts are monotonic (and increasing) on $I=[\mu_1,\mu_2]$ if 
\begin{equation}\label{eq:vk}
v_k= [\min_I \phi_k^\prime (\mu)]-\max_{j\in\{k-1,k\}}\max_I s_j(\mu )\, [\max_I b^\prime_k(\mu)]>0, \quad k=1, \dots ,n.
\end{equation}

\begin{lemma}\label{lem:linint}
Suppose that the lifts $\{F_\mu\}$ of a family of PWL circle maps are defined by $(b_k(\mu ), \phi_k(\mu ))$, $k=1,\dots ,n$, for $\mu\in I=(\mu_1,\mu_2)$. Suppose further that for all $\mu,\nu\in I$, $b_k(\mu )\in (b_{k-1}(\nu ), b_{k+1}(\nu)$ and $b_k(\nu )\in (b_{k-1}(\mu ), b_{k+1}(\mu))$. If $b_k$ and $\phi_k$ are $C^2$ functions of $\mu$ and  
\begin{equation}\label{eq:minmaxcond}
C_1= \min_k v_k >0, 
\end{equation}
where $v_k$ is defined in (\ref{eq:vk}) then $\{F_\mu\}$ is the lift of a monotonic family of circle maps and for any  $\mu >\nu$ in $I$ there exists $C_2>0$ such that for all $x\in{\mathbb{R}}$
\begin{equation}\label{eq:linlowbound}
C_1(\mu -\nu)<F_\mu (x)- F_\nu (x)< C_2(\mu-\nu).
\end{equation}
\end{lemma}

\emph{Proof:} Since $F_\mu$ is PWL the minima and maxima of $F_\mu (x)-F_\nu (x)$ occur at break points of the maps. Note that (\ref{eq:Fk}) implies that the branch of $F_\mu$ in $(b_{k-1},b_k)$ can also be written as
\begin{equation}\label{eq:Fk-1}
F_\mu (x)=\phi_k(\mu )+s_{k-1}(\mu)(x-b_k(\mu )), \quad x\in [b_{k-1}(\mu),b_k(\mu )].
\end{equation}
The monotonicity condition $F_\mu (x)-F_\nu (x)>0$ for all $x$ holds provided it holds at each break point of $F_\mu$ and at each break point of $F_\nu$.   

Suppose that $b_k(\nu )<b_k(\mu )$. 

Then $F_\nu (b_k(\nu ))=\phi_k (\nu)$ by definition and using (\ref{eq:Fk-1})
\[
F_\mu (b_k(\nu ))=\phi_k (\mu)+s_{k-1}(\mu )(b_k(\mu )- b_k(\nu)).
\]
Hence
\[
F_\mu (b_k(\nu ))-F_\nu (b_k(\nu ))=\phi_k (\mu)-\phi_k(\nu )+s_{k-1}(\mu )(b_k(\nu )- b_k(\mu))
\]
and so by the mean value theorem
\begin{equation}\label{eq:mvt}
F_\mu (b_k(\nu ))-F_\nu (b_k(\nu ))=\phi_k' (r_1)(\mu-\nu )-s_{k-1}(\mu )b_k'(r_2)(\mu )-\nu)),
\end{equation}
for some $r_i\in [\nu,\mu ]$, $i=1,2$. By (\ref{eq:minmaxcond}) the right hand side of (\ref{eq:mvt}) is greater than or equal to $C_1(\mu -\nu)$ as required. The case of $F_\mu (b_k(\mu ))-F_\nu (b_k(\mu ))$ is similar, as are the remaining two cases for $b_k(\mu )<b_k(\nu )$.

The upper bound is proved by a small modification of this argument.
\newline\rightline{$\square$}

Equation (\ref{eq:minmaxcond}) of Lemma~\ref{lem:linint} is a global condition. The proof implies the following local corollary.

\begin{corollary}\label{cor:linint}
Suppose that $F_\mu$ are the lifts of a family of PWL circle maps $\{f_\mu\}$ defined by $(\phi_k,b_k)_1^n$ and
\begin{equation}\label{eqlem:linint}
\phi'_k(\mu_c)-\max\{0,s_{k-1}(\mu_c)b_k'(\mu_c),s_{k-1}(\mu_c)b_k'(\mu_c)\} >0, \quad k=1, \dots ,n,
\end{equation}
then there exists an open neighbourhood of $\mu_c$ in which $\{f_\mu\}$ is a monotonic  family.
\end{corollary}
 
Note that in both results, Lemma~\ref{lem:linint} and Corollary~\ref{cor:linint}, the conditions simplify if the break points are fixed to simply having $\phi_k'(\mu)>0$.

To describe properties of families of maps close to parameters at which there is a conjugacy to a rational rigid rotation with rotation number $\frac{p}{q}$, it is useful to understand the structure of $F_\mu^q$ near that parameter. Clearly at the critical parameter, $\mu =0$ say, $F_0^q(x)=x+p$, and without loss of generality we may work with a lift $F^q$ having $p=0$ provided the trivial integer is added into the result for $F_\mu$ at the end.  

Break points of $f_\mu^q$ are contained in the set of preimages of the break points of $f_\mu$, or more precisely the are in the set
\begin{equation}\label{eq:betakj}
B(\mu )=\{\beta_{kj}(\mu )~|~f_\mu^{-j}(\beta_{kj})=b_k (\mu), \ 0\le j<q,\ 1\le k\le n\}.
\end{equation}
At a parameter value at which there is a conjugacy to a rigid rational rotation $f^q$ has no break points (it equals the identity) but the set $B$ is the union of the orbits of the break points. Since at least two break points are on each orbit, as a parameter varies we can expect the coincidences (a point being the preimage of at least two break points) to become disentangled. The next couple of results describe this process, and show that the break points remain in order $\mu$ neighbourhoods of their position at $\mu =0$.  

Let $N(x,\delta )$ denote the open all of radius $\delta$ about $x$.

\begin{lemma}\label{lem:preims}Suppose that $\{f_\mu \}$ is a 
monotonic continuous family of orientation preserving circle homeomorphisms and $f_0$ is conjugate to a rigid rational rotation $x\to x+\frac{p}{q}$. If $b_k(\mu )$ and $\phi_k(\mu )$ are $C^2$ functions then there exists an open neighbourhood ${\mathcal N}$ of $\mu=0$ and $\gamma >0$ such that
\[
B(\mu )\subset \cup_{j,k}N(\beta_{j,k}(0),\gamma |\mu |).
\]
\end{lemma}

\emph{Proof:} If $\mu =0$ then each break point of $f_0$ is on the orbit of at least one other break point and so the points $\beta_{jk}(0)$ are not distinct. Applying Lemma~\ref{lem:linint} to $F^{-1}$ there is a neighbourhood ${\mathcal N}$ of $\mu =0$ such that
$|F_\mu^{-1}(x)-F_0^{-1}(x)|<K_1|\mu|  $ for all $x\in {\mathbb{R}}$ and $\mu\in {\mathcal N}$. Moreover, $|F_\mu^{-1} (y)-F_\mu ^{-1}(x)|< K_2|y-x|$ where $K_2=\max\{1/s_k(\mu)\}$. Hence $|F_\mu^{-2}(x)-F_0^{-2}(x)|<K_1K_2|\mu |$ and by induction  $|F_\mu^{-r}(x)-F_0^{-r}(x)|<K_1K_2^{r-1}|\mu |$ and the result follows.
\newline\rightline{$\square$}

\section{Scaling of the rotation number}\label{sect:scaling}

\medskip
\emph{Proof of Theorem~\ref{thm:scaling}:} Without loss of generality assume that $\mu_c=0$ and the family $F_\mu$ is increasing. By assumption there is a lift $F_0$ of $f_0$ such that $\rho (F_0)=\frac{p}{q}$, $f_0^q(x)=x$ and $F^q_0(x)=x+p$. Since lifts differ by an integer we may choose a family of lifts  
\begin{equation}
G_\mu (x)=F^q_\mu (x)-p
\label{eq:G}
\end{equation}
of $f_\mu^q$ such that $G_0(x)=x$. Note that $G_\mu$ is not the $q^{th}$ iterate of a lift of $f_\mu$ (the $q^{th}$ iterates of lifts of $f_\mu$ differ by an integer multiple of $q$) , but $\rho (G_0)=0$ and $\rho (F^q_\mu)=p+\rho (G_\mu)$, so 
\begin{equation}
\rho (F_\mu )=\frac{p}{q}+\frac{1}{q}\rho (G_\mu ).
\label{eq:rhoG}
\end{equation}

Now, $|G_\mu(x)-G_0(x)|=|G_\mu(x)-x|$ and so by the proof of Lemma~\ref{lem:preims} there is a neighbourhood of $\mu =0$ and constants $C_3>0$ and $C_4>0$ such that 
\begin{equation}\label{eq:sc1}
C_3\mu <G_\mu (x)-x<C_4 \mu 
\end{equation} 
for all $\mu$ in this neighbourhood. This implies that $|\rho(\mu )-\frac{p}{q}|>C_5\mu$
for some $C_5>0$ immediately.

\medskip
The proof of (\ref{eq:sc2}) requires a little more work. We start by restricting the preimages of the break points of $f$ which can be break points of $f^q$ in size $|\mu |$ neighbourhoods $L_i$ of the orbits of the break points of $f$ at $\mu =0$. Outside these neighbourhoods $f^q$ has slope close to one and evolves smoothly. If $\mu >0$ then an estimate of the time spent by an orbit outside the neighbourhoods $L_i$ can be made (Lemma~\ref{lem:timei}) and the time in each $L_i$ is bounded (Lemma~\ref{lem:btimei}) and so does not contribute to the rotation number. These estimates can be used to show that the rotation number scales linearly in $\mu >0$ and finally a symmetry argument shows that the scaling is the same in $\mu <0$, completing the proof.  

Recall that $K$ is the number of distinct periodic orbits containing the break points at $\mu =0$ (cf. Theorem~\ref{thm:Jcon}). Thus by choosing a smaller neighbourhood of $\mu =0$ the neighbourhoods $N(\beta_{jk}(0), \gamma |\mu|)$ of Lemma~\ref{lem:preims} are disjoint. There are clearly $qK$ such intervals, so label these $L_i$, $1\le i\le qK$ monotonically around the circle so that $x\in L_i$, $y\in L_{i+1}$ implies that $x<y$ with the obvious modification if $i=qK$. If $\beta_{kr}(0)\in L_i$ define
\begin{equation}\label{defell}
\ell_i=\beta_{kj}(0),
\end{equation}
and so
\begin{equation}\label{eq:defmiMi}
m_i(\mu)=\ell_i-\gamma |\mu|, \quad M_i(\mu)=\ell_i+\gamma |\mu|.
\end{equation}

Between $M_i$ and $m_{i+1}$ there are, by definition, no break points of $f^q$ and hence $f^q$ is close to the identity and for any $y_i\in (M_i,m_{i+1})$ the lift $G_\mu$ of (\ref{eq:G}) is defined by
\begin{equation}
G_\mu (x)=G_\mu (y_i)+S_i(x-y_i), \quad x\in [M_i,m_{i+1}],
\label{eq:linmap}
\end{equation}
where $S_i(\mu )$ is the product of $q$ of the piecewise linear slopes $s_k$ of (\ref{eq:sk}) and $S_i(0)=1$. Hence $|S_i-1|<C_6 |\mu |$ (since there are only a finite number of these segments, $C_6$ can be chosen to be independent of $i$). The bound of $|S_i-1|$ as a function of $|\mu |$ follows from the $C^2$ variation of $(\phi_k)$ and $\beta_k$ together with (\ref{eq:sk}). 

For sufficiently small $\mu_0$ we may choose $y_i$ so that it is contained in $(M_i(\mu_0),m_{i+1}(\mu_0))\cup (M_i(-\mu_0),m_{i+1}(-\mu_0))$ and hence in $(M_i,m_{i+1})$ for all $|\mu |<\mu_0$. Since $f$ is piecewise $C^2$ this implies that $G_\mu (y_i)$ is $C^2$ in $(-\mu_0,\mu_0)$ and in particular that 
\begin{equation}
\frac{\partial G_\mu}{\partial \mu}(y_i)\big|_{\mu=0} \ {\rm exists\ and \ equals}\ A_i>0, \ {\rm say}.
\label{eq:partialFq}
\end{equation}
The increasing property, $A_i>0$ follows from the increasing property of $F$.  

\begin{figure}[t!]
\centering
\includegraphics[width=10cm]{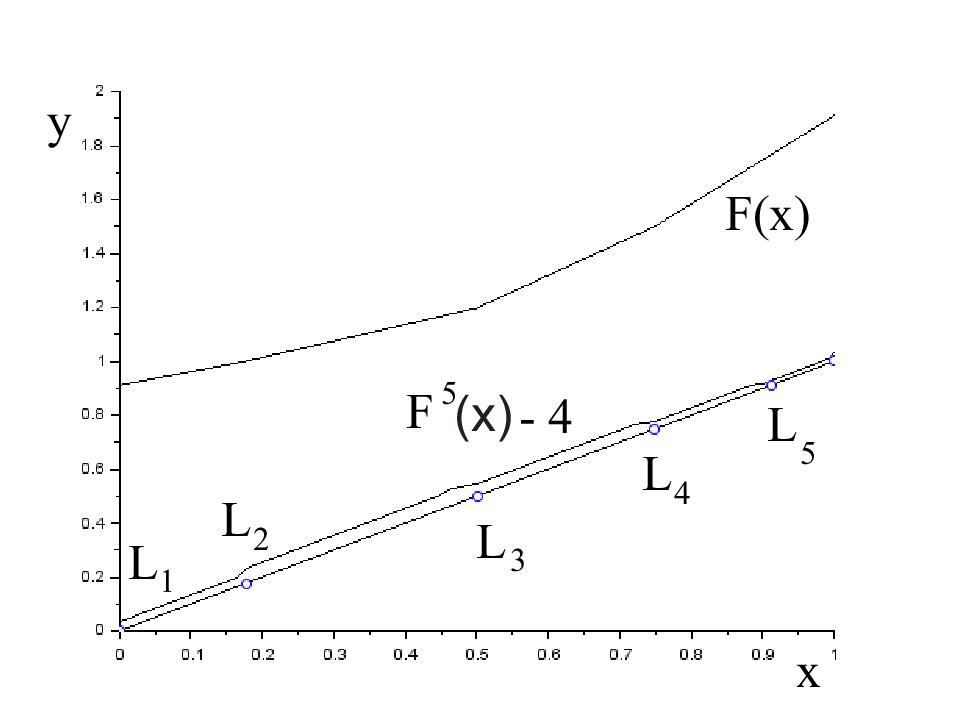}
\caption{
Sketch of a lift $f_{\alpha ,\beta}$  and the iterate $f_{\alpha ,\beta}^5$ shifted by 4 using (\ref{eq:GMMex}) with $\alpha =2$ and $\beta=(\sqrt{5}-1)-\mu$, $\mu=0.02$ (cf. Figure~\ref{fig:ourfig}a). As discussed in section~\ref{sect:GMMex}, $f_{\alpha ,\beta}$ has four break points and if $\mu =0$ these lie on the same orbit of period five which is indicated by the circles on the diagonal. The intervals $L_i$ are centred on these (the point which is not a break point is in $L_4$, and contain all the break points of $F_{\alpha ,\beta}^5$. Outside these regions the map has slope equal to one (in general the slope would be close to one) Note that $L_1$ is centred on $x=0$ and so straddles the periodic window presented here.  
}
\label{fig:msketches}
\end{figure}

Thus if $|\mu |$ is small the map has the features shown in Figure~\ref{fig:msketches}: over most of the circle the map is linear with slope close to 1, but in the $qK$ order $\mu$ intervals $L_i$ the map can have different slopes not necessarily close to 1. 

Now assume that $\mu>0$, so $G_\mu (M_i)>C_7\mu$ for some $C_7>0$ by the proof of Lemma~\ref{lem:preims}. 

To estimate the rotation number of $f^q$ (and hence $f$) the length of time spent in these order one intervals can be added together, and the $qK$ order $\mu$ intervals between them will contribute a bounded extra contribution which will not add to the rotation number.

\medskip
\begin{lemma}\label{lem:timei}Let $G_{\mu,i}=G_\mu (M_i)$ and $\mu >0$. If $x\in [M_i,G_{\mu ,i}]$ then the number of iterates of $x$ in  $[M_i,m_{i+1}]$ is $T_i\in [n_i-2,n_i]$ where
\begin{equation}\label{eq:nidef}
G_{\mu,i}+S_i(1+S_i+\dots +S_i^{n_i-1})(G_{\mu,i}-M_i) = m_{i+1}.
\end{equation}
Moreover, $G_{\mu,i}-M_i=A_i\mu +O(\mu^2)$ and $n_i=\frac{\kappa_i}{\mu}+O(1)$ where if $S_i=1+O(\mu^2)$ and $A_i$ and $\kappa_i$ are constants with
\begin{equation}\label{eq:n1s1}
\kappa_i=\frac{\ell_{i+1}-\ell_i}{A_i},
\end{equation}
and if $S_i=1+B_i\mu +O(\mu^2)$ then 
\begin{equation}\label{eq:n1s2}
\kappa_i=\frac{1}{B_i}\ell n\,\left(1+\frac{(\ell_{i+1}-\ell_i)B_i}{A_i}\right) .
\end{equation}
\end{lemma}

\emph{Proof:} Choose $y_i$ as described above equation (\ref{eq:partialFq}). Then
$G_\mu (M_i)-M_i=G_{\mu,i}-M_i=G_\mu (y_i)-S_iy_i + (S_i-1)M_i$. By assumption $|S_i-1|<C_6\mu$ and observe that $G_\mu (y_i)-S_iy_i$ is independent of $y_i$ (it is $G_\mu (0)$ for the continuation of the function defined by (\ref{eq:linmap})). This is therefore bounded above and below by multiples of $\mu$ by the same argument as in the proof of Lemma~\ref{lem:preims}, so $G_{\mu,i}-M_i=A_i\mu +O(\mu^2)$, $A_i>0$. 

Iterating $G_\mu$, shifted to start at $M_i$, $n$ times gives 
\[
G_\mu^{k}(x,\mu )=M_i+(G_{\mu,i}-M_i)(1+S_i +\dots +S_i^{k-1})+S_i^{k}(x-M_i).
\]
Let $k$ be the smallest integer such that $G_\mu^{k}(M_i,\mu )\le m_{i+1}<G_\mu^{k+1}(M_i,\mu )$. Such a $k$ exists by (\ref{eq:sc1}) and in particular there are no fixed points of $G_\mu$ in $[M_i,m_{i+1}]$. This means that the quantities below make sense and in particular the argument of the logarithms is never negative.  Thus hen the orbit of $M_i$ has $k$ points in $[M_i,m_{i+1}]$ and $G_{\mu ,i}$ has $k-1$. Hence for all $x\in [M_i,G_{\mu ,i}]$ there are either $k$ or $k-1$ such points.

Let $n_i$ be the solution to 
\begin{equation}\label{eq:Fni}
(G_{\mu,i}-M_i)(1+S_i +\dots +S_i^{n_i-1})=m_{i+1}-M_i.
\end{equation}.
Then $n_i-1<k\le n_i$ and so the number of iterates of a point $x\in [M_i,G_{\mu ,i}]$ which lie in in $[M_i,m_{i+1}]$ is in $[n_i-2,n_i]$.   

If $S_i=1+O(\mu^2)$ then the solution of (\ref{eq:Fni}) is 
\begin{equation}\label{eq:nis=1}
n_i=\frac{m_{i+1}-M_i}{G_{\mu ,i}-M_i}=\frac{\ell_{i+1}-\ell_i}{A_i\mu}+O(1).
\end{equation}
By assumption $G_{\mu ,i}-M_i=A_i\mu+O(\mu^2)$, so the definition of $\kappa_i$ (in the statement of Lemma~\ref{lem:timei}: $n_i=\frac{\kappa_i}{\mu}+O(1)$) implies (\ref{eq:n1s1}).

If $S_i=1+B_i\mu+O(\mu^2)$ then the left hand side of (\ref{eq:Fni}) is 
\[
(G_{\mu ,i}-M_i)\frac{S_i^{n_i}-1}{S_i-1}=\frac{A_i}{B_i}((1+B_i\mu +O(\mu^2))^{n_i})-1) +O(\mu ).
\]
and since $n_i=\frac{\kappa_i}{\mu}+O(1)$,  
\begin{equation}\begin{array}{rl}
(1+B_i\mu +O(\mu^2))^{n_i} &=\left(1+B_i\mu+O(\mu^{2})\right)^{\frac{\kappa_i}{\mu}}\\
 &= e^{B_i\kappa_i}\left(1+O(\mu )\right)
\end{array}\end{equation}
as $\mu \to 0$. Putting these together (\ref{eq:Fni}) becomes to lowest order
\[A_i (e^{B_i\kappa_i}-1)=B_i(m_{i+1}-M_i )+O(\mu),\]
and so a little algebra along with the relation $m_{i+1}-M_i=\ell_{i+1}-\ell_i+O(\mu )$ gives (\ref{eq:n1s2}).
\newline\rightline{$\square$}

Note that the limit $B_i\to 0$ of (\ref{eq:n1s2}) is (\ref{eq:n1s1}) as expected.

Lemma~\ref{lem:timei} describes the time (number of iterates) taken on the `long' passages where $G_\mu $ has slope close to one. The remaining time is spent in the order $\mu$ regions $L_i$ of Figure~\ref{fig:msketches} close to the orbits of the break points at $\mu=0$, i.e. between $m_i$ and $M_i$.  

\begin{lemma}There exists $D>0$ independent of $\mu$ such that the number of iterates $(G_\mu^k(0))$, $G_\mu^k(0)<1$, in $\cup [m_i,M_i]$ is bounded by $D$ for all $\mu\in{\mathcal N}$.
\label{lem:btimei}
\end{lemma}

\emph{Proof:} The inequality (\ref{eq:sc1}) $C_3\mu <G_\mu (x)-x<C_4 \mu$  implies that $G_\mu^{k}(x)\ge x+kC_3\mu$. Thus at most $k_i+1$ iterates can lie in $[m_i,M_i]$ on each passage through the circle where $k_i=(M_i-m_i)/(C_3\mu )$. Since $M_i-m_i=2\gamma\mu$ the total length is bounded by $2qK\gamma\mu$ (where $K$ is the number of different orbits of the break points at $\mu =0$) and the number of iterates taken is bounded by the sum of the $k_i$: $D=(2qK\gamma /C_3)+qK$ which is independent of $\mu$.  
\newline\rightline{$\square$}

\medskip
\emph{Completion of the proof of Theorem~\ref{thm:scaling}}

Suppose that $\mu >0$. Let $N$ be the number of iterates in a complete passage through a distance one. Then
\[
\sum (n_i-2)< N < D+\sum n_i.
\]
Hence, using (\ref{eq:rhoG}) to reinstating the integer $p$ into the lift of $F_\mu^q$,
\begin{equation}\label{eq:rotqF}
p+\frac{1}{D+\sum n_i}\le \rho (F_\mu ^q)\le p+\frac{1}{\sum (n_i-2)}.
\end{equation}
Now, in both cases the denominators are dominated by $\sum n_i$. Indeed, since $n_i=\frac{\kappa_i}{\mu}+O(1)$,
\[
\frac{1}{D+\sum n_i}=\frac{\mu}{\sum \kappa_i}\left(1+O(\mu)\right), \quad \frac{1}{\sum (n_i-2)}=\frac{\mu}{\sum \kappa_i}\left(1+O(\mu)\right).
\]
Hence
\[
\rho (F_\mu )-\frac{p}{q}-\frac{\mu}{q\sum \kappa_i} =O(\mu^2).
\]
This is (\ref{eq:sc2}) of Theorem~\ref{thm:scaling} with $R_1$ given by (\ref{eq:R1}).

We now need to establish the same result \emph{with the same linear growth condition} in $\mu<0$. This follows from a symmetry between the two cases, or more precisely an `almost' symmetry: in the long laminar regions between $M_{i}$ and $m_{i+1}$, $G_\mu$ is approximately $G^{-1}_{-\mu}$. This implies that all the scaling results for $\mu >0$ hold at $-\mu$ but with the direction of passage reversed.

To prove this first note that in these laminar regions the map is a composition of affine maps whose coefficients depend $C^2$ on the parameter $\mu$, so that is also true of $F^q$. Indeed, in these regions (\ref{eq:linmap}) and (\ref{eq:partialFq}) imply 
\[ G_\mu (x )=y_i+A_i\mu+S_i(x-y_i),\]
valid in a neighbourhood of $\mu=0$. So $S_i=1+B_i\mu+O(\mu^2)$ implies that
\begin{equation}
G_\mu (x )=x+\mu (A_i+B_i(x-y_i)+O(\mu^2).
\label{eq:mudepFq}
\end{equation}
Equation (\ref{eq:mudepFq}) is of the form $G(x)=x+\mu g(x)+O(\mu^2)$, for which, using the method of successive approximation,  $G^{-1}(x)=x-\mu g(x)+O(\mu^2)$, i.e. if $\mu >0$ then
\begin{equation}
G_{-\mu}(x )=G_\mu^{-1}(x)+O(\mu^2).
\label{eq:symm}
\end{equation}
This implies that to lowest order the estimates of Lemma~\ref{lem:timei} and Lemma~\ref{lem:btimei} hold in $\mu<0$ \emph{with the same coefficients} but with the orbit travelling in the opposite direction.
(Or equivalently, the rotation number of the inverse is minus the rotation number of the function.) Hence the symmetry up to the terms of order $\mu$ (\ref{eq:symm})  implies that $\rho (G_{-\mu})=-\rho (G_\mu )+O(\mu^2)$ and the symmetry of the constant $R_1$ of (\ref{eq:sc2}) follows.

This completes the proof of Theorem~\ref{thm:scaling}.
\newline\rightline{$\square$}
Note that the strategy of the argument, separating the analysis of the orbit into small regions containing the break points and larger intervals between them, is similar to the construction used on a piecewise smooth example by Quas \cite{Quas1992}. 

\section{Herman's Example}\label{sect:Hex}
Herman \cite{H1979} discusses the PWL map with lift defined $[0,1)$ by
\begin{equation}\label{ex:H}
x_{n+1}=H(x_n) =\mu +h(x_n), \quad h(x)= \begin{cases}\lambda x & {\rm if}~0\le x\le c\\ 1+\lambda^{-\beta} (x-1) & {\rm if}~c<x<1\end{cases}
\end{equation}
and extended to the $\mathbb{R}$ using $H(x+1)=H(x)+1$. The parameters are chosen so that
\begin{equation}\label{eq:Hconds}
\lambda >1, \quad \beta >0, \quad \lambda c=1+\lambda^\beta (c-1).
\end{equation}

Herman \cite{H1979} (see also \cite{Dzh2018}) proves the following theorem.

\begin{theorem}\label{thm:H1979}Suppose $H$ has irrational rotation number $\rho$. Then the following are equivalent:
\begin{itemize}
\item[(i)] The break points $0$ and $c$ lie on the same orbit;
\item[(ii)] $H$ is Lipschitz conjugate to a rigid irrational rotation;
\item[(iii)]  $\frac{\beta}{1+\beta}\in \rho\mathbb{Z}$.
\end{itemize}\end{theorem}

Coelho et al \cite{Coelho1995} consider a sub-family of Herman's example which satisfies condition (i) automatically and allows us to consider the case of rational rotation numbers for Herman's example. For $a$ and $b$ in $(0,1)$ consider 
\begin{equation}\label{eq:defCoelho}
G(x)=\begin{cases}a+\frac{1-a}{b}x & {\rm if}~0\le x<b,\\
1+\frac{a}{1-b}(x-b) &{\rm if}~b\le x<1.
\end{cases}
\end{equation}
For this family there is no mode locking \cite{Coelho1995}, which also follows from our Theorem~\ref{thm:nomodelock}. The key to this example is that $G(b)=1$ so the break points $0$ and $b$ lie on the same orbit of the associated circle map (and hence are conjugate to a rigid rational rotation if the rotation number is rational). If the slopes of the map are
\[
\alpha_G=   \frac{1-a}{b}, \quad \beta_G=\frac{a}{1-b},
\]
then \cite{Coelho1995}
\begin{equation}\label{eq:rotG}
\rho (G)=\frac{\ell n\, \alpha_G}{\ell n\, \alpha_G-\ell n\, \beta_G}.
\end{equation}

Theorem~\ref{thm:H1979} and (\ref{eq:rotG}) make it possible to construct families of maps passing through maps with rational rotation number conjugate to a rigid rotation. 

The example we have chosen to illustrate Theorem~\ref{thm:scaling} is given below, note that the parameter $\mu$ of (\ref{ex:H}) has been shifted by a constant $\frac{1}{1+\lambda}$ so that $\mu =0$ will correspond to a member of the family (\ref{eq:defCoelho}).

\begin{figure}[t!]
\centering
\includegraphics[height=7.5cm]{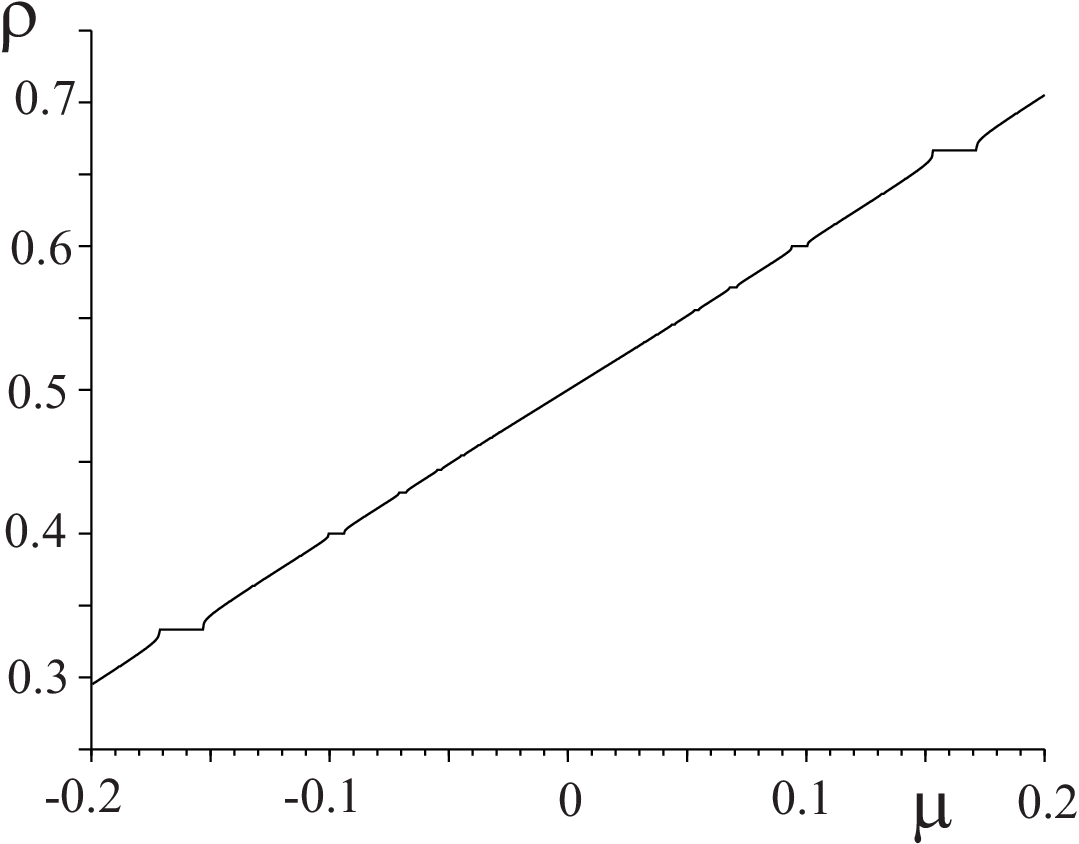}
\caption{
Rotation number of (\ref{def:R}) as a function of the parameter $\mu$. The rotation number is calculated using 100000 iterates of the map at 1000 equidistant values of the parameter between $-0.2$ and $0.2$.
}
\label{fig:Hfig1}
\end{figure} 

\begin{equation}\label{def:R}
x_{n+1}=R(x_n)=\mu +\frac{1}{1+\lambda} + Y(x_n), \quad Y(x)= \begin{cases}\lambda x & {\rm if}~0\le x\le c\\ 1+\lambda^{-1} (x-1) & {\rm if}~c<x<1\end{cases},
\end{equation}
with $c=\frac{1}{1+\lambda}$ to satisfy (\ref{eq:Hconds}). This corresponds to the case $\beta=1$ of (\ref{ex:H}). If $\mu =0$ this is also one of the Coelho et al \cite{Coelho1995} maps, (\ref{eq:defCoelho}), with $\alpha_G=\lambda$ and $\beta_G=\lambda^{-1}$ (i.e. $a=b=\frac{1}{1+\lambda}$) and so $\rho (0) =\frac{1}{2}$. Figure~\ref{fig:Hfig1} shows the rotation number as a function of $\mu$ with $\lambda=\sqrt{2}$. Here and in section~\ref{sect:GMMex} the computations are done in the {\tt Scilab} environment \cite{Scilab}. It certainly looks like a small perturbation from a straight line (some mode locking regions with the standard logarithmic scaling are visible) and linear regression using the {\tt Scilab} \emph{reglin} command, shows that the best linear fit is $\rho = 0.500+1.027\mu$ to three decimal places, with standard deviation $0.002$. The constant $0.500$ is the rotation number at $\mu=0$ and the slope can be calculated from (\ref{eq:R1}). Clearly $f(\mu,0)=f(\mu,c)=1$ and so the change in the second iterate is the product of these, $A_i= 1$, and the second iterate has no other break points if $\mu =0$. Hence by (\ref{eq:R1}) $\sum \kappa =(\ell_0+1-\ell_1)+(\ell_1-\ell_0)=1$ and the slope is unity, as confirmed by the numerical simulations.

\section{An example from refraction in a periodic medium}\label{sect:GMMex}
In this section we consider an example described in \cite{GMM2024}. As a function of the circle it can be written as
\begin{equation}\label{eq:GMMex}
f_{\alpha ,\beta}(x)=\begin{cases}
\frac{1}{\alpha}x+\frac{1}{2\beta}(\beta+1)&{\rm if} ~0\le x<\frac{\alpha}{2\beta}(\beta -1)\\
\frac{\beta}{\alpha}x+\frac{1}{2}(1-\beta )&{\rm if} ~\frac{\alpha}{2\beta}(\beta -1) \le x<\frac{1}{2}\\
\beta(x-1)+\frac{1}{2\alpha}(\alpha +\beta )&{\rm if} ~\frac{1}{2}\le x< 1-\frac{1}{2\alpha}\\
\frac{\alpha}{\beta}(x-1)+\frac{1}{2\beta}(\beta +1)&{\rm if} ~1-\frac{1}{2\alpha}\le x <1\end{cases}.
\end{equation}
The parameters are restricted to 
\begin{equation}\label{eq:gmmconds}
\alpha >1, \quad \beta >1, \quad \alpha^{-1}+\beta^{-1} > 1.
\end{equation}
This system contains a number of symmetries which make it possible to characterise the dynamics through renormalization in a very nice way. Here we will only consider the local scaling near a parameter value with conjugacy to a rigid rational rotation number; see \cite{GMM2024} for global results and a finer description of the scaling near these points.

Although there are four break points,
\[
b_1=0, \quad b_2= \frac{\alpha}{2\beta}(\beta -1), \quad b_3=\frac{1}{2} , \quad {\rm and} \quad b_4=1-\frac{1}{2\alpha},
\]
their orbits are not independent. Direct calculation shows that 
\[
f_{\alpha ,\beta}(b_2)=b_1+1=b_1~({\rm mod} ~1), \quad f_{\alpha ,\beta}(b_4)=b_3.
\]
Example (\ref{eq:GMMex}) has symmetries which imply that conjugacy to rigid rational rotations becomes a codimension one phenomenon. In particular, for the case of Theorem~\ref{thm:gmmrot} below, $f_{\alpha,\beta}(b_3)=b_2$ if and only if $f^2_{\alpha ,\beta }(b_1)= b_4$ and hence only one condition is needed to ensure the existence of a periodic orbit of period five and rotation number $\frac{4}{5}$ involving all the break points and one other point, $f_{\alpha ,\beta}(b_1)$.

Note that with the parametrization of (\ref{eq:GMMex}) the families considered are increasing continuous families as $\mu$ decreases, hence the rotation number is a decreasing function of $\mu$. 

\begin{theorem}\label{thm:gmmrot}If $\alpha >1$ and $\beta =\frac{--\alpha+\sqrt{\alpha^2(\alpha-1)^{-1}(\alpha+3)}}{2}$ then (\ref{eq:GMMex}) has a periodic orbit of period five which contains all four break points. At these parameters the map is conjugate to a rigid rotation with rotation number $\frac{4}{5}$.\end{theorem}

\emph{Proof:} By direct calculation $f_{\alpha, \beta}(b_3)=b_2$ if $(\alpha -1)\beta^2+\alpha (\alpha -1)\beta -\alpha^2=0$ which has a positive solution for $\beta$ with $\beta=\frac{\sqrt{5}-1}{2}\alpha$. The condition $\alpha^{-1}+\beta{-1}>1$ of (\ref{eq:gmmconds}) implies that $ f_{\alpha ,\beta}(1)>b_4$ and hence that $f^2_{\alpha,\beta}(1)=\frac{1}{2\beta^2}(\alpha +\beta -\alpha \beta +\beta^2)$ and setting  $f^2_{\alpha,\beta}(1)=b_4$ gives the same quadratic equation for $\beta$ as the condition $f_{\alpha, \beta}(b_3)=b_2$. The solution to the quadratic is 
\[
\beta =  \frac{-\alpha+\sqrt{\alpha^2(\alpha-1)^{-1}(\alpha+3)}}{2},
\]
and two elementary calculation show that the second two conditions of (\ref{eq:gmmconds}) are automatically satisfied if $\alpha >1$.

Thus the points $0$, $f_{\alpha,\beta}(0)$, $b_4$, $b_3$, and $b_2$ map to each other forming an orbit of period five with rotation number $4/5$, and hence the theorem follows from Theorem~\ref{thm:equ}.
\newline\rightline{$\square$}

\begin{figure}[t!]
\centering
\includegraphics[height=5.5cm]{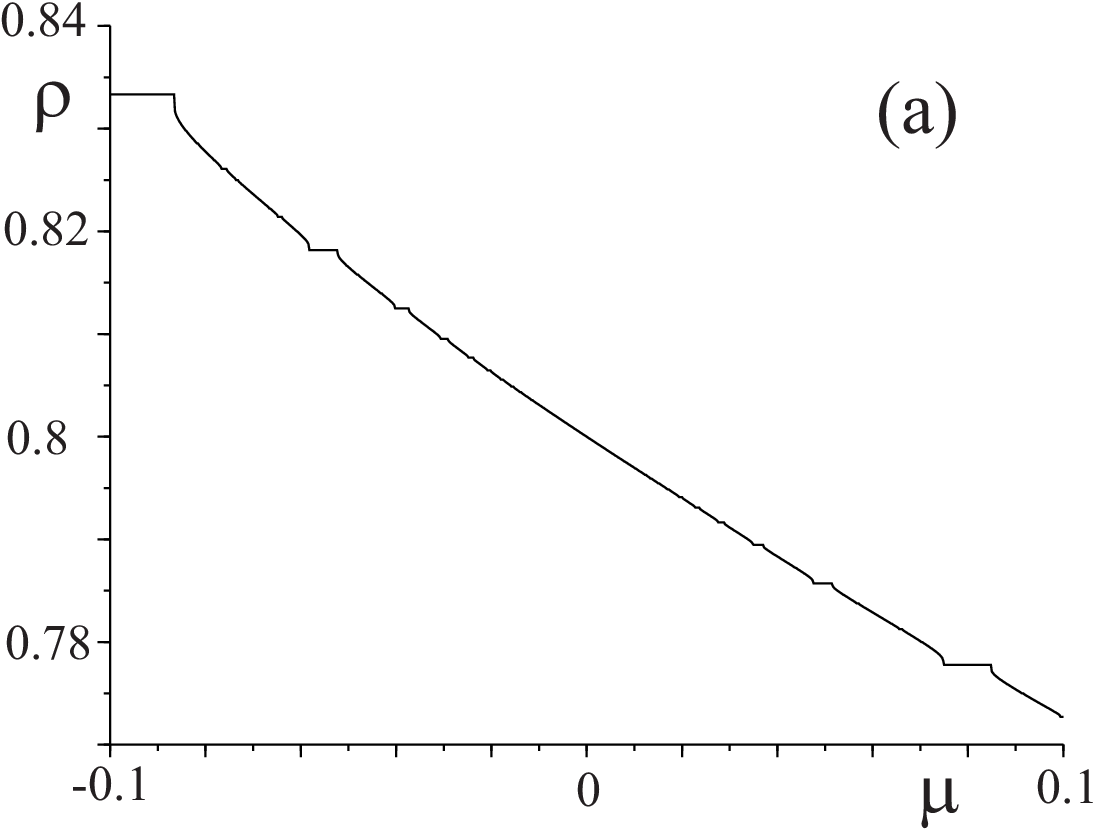}
\includegraphics[height=5.5cm]{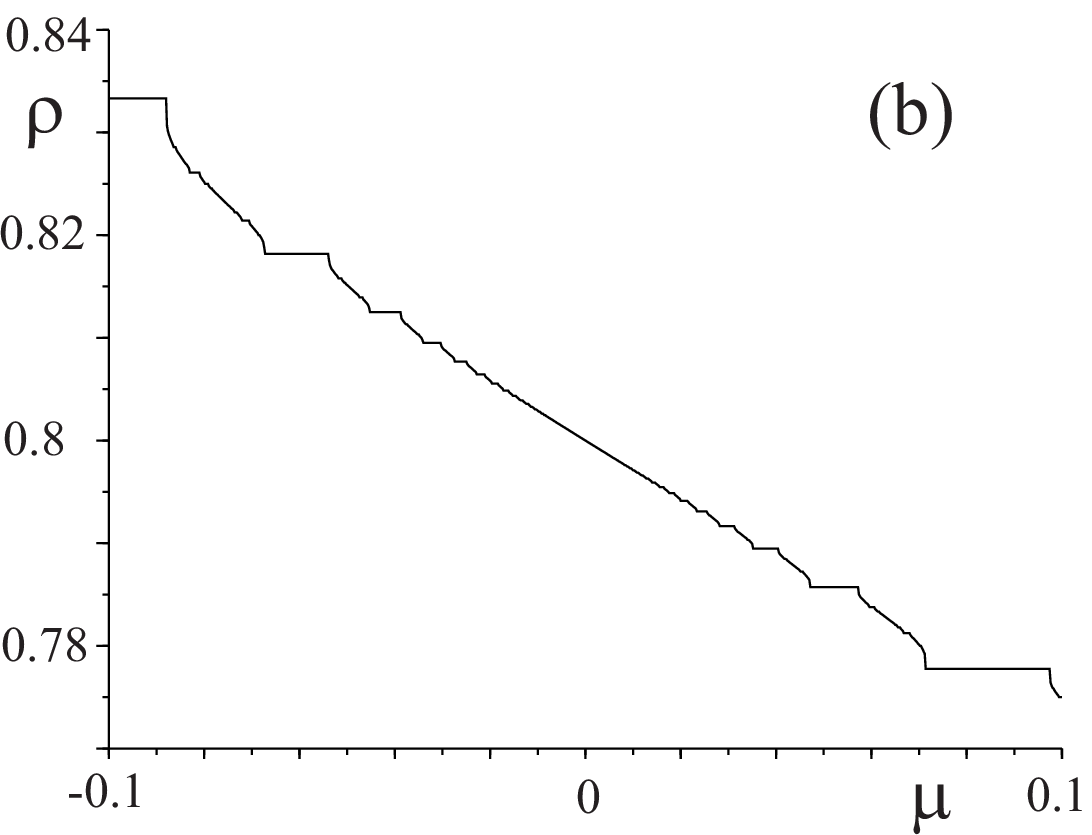}
\caption{
Rotation number (vertical axis) of (\ref{eq:GMMex}) against the parameter $\mu$. The rotation number is calculated using 100000 iterates of the map at 1000 equidistant values of the parameter $\mu$. (a) $\alpha=2$, $\beta=(\sqrt{5}-1)+\mu$, $-0.1<\mu <0.1$; (b) $m=\frac{\alpha}{\beta}=\sqrt{10}$, $\alpha =\alpha_0+\mu$, $\alpha_0= 3.403$, $-0.1<\mu<0.1$. 
}
\label{fig:ourfig}
\end{figure} 

To demonstrate the scaling of Theorem~\ref{thm:scaling} in Figure~\ref{fig:ourfig}a we take $\alpha =2$ and calculate $\rho (f_{\alpha ,\beta})$ with $\beta$  varying through $\beta_0=\sqrt{5}-1\approx 1.236$. The scaling is linear and linear
regression suggests the linear approximation is $\rho \sim 0.801-0.312(\beta -\beta_0)$ with standard deviation $0.0014$. By choosing a smaller range of perturbations the rotation number at $\beta=\beta_0$ from interpolation can be brought much closer to $4/5$ and the standard deviation reduced as well. Because of the greater complexity of this map we have not compared the slope of the graph with the slope predicted by (\ref{eq:R1}). Note that there is an obvious global variation of the graph from linear.  In \cite{GMM2024} a different parameterisation is used: the ratio $m=\frac{\alpha}{\beta}$ is held fixed. With the choice $m=\sqrt{10}$ of \cite{GMM2024}, Theorem~\ref{thm:gmmrot} shows that the system is conjugate to a rigid rational rotation with rotation number $4/5$ if $\alpha =\alpha_0\approx 3.403$ (with corresponding $\beta =\alpha_0/m\approx 1.076$. The numerically calculated rotation number as a function of $\alpha$ around this value, with $m$ fixed, is shown in Figure~\ref{fig:ourfig}b. Again, linear regression shows $\rho =0.801-0.294(\alpha -\alpha_0)$ with standard deviation $0.002$. The standard mode-locking plateaus are particularly clear in this diagram, although there are further values of parameters with systems conjugate to rigid rotations within this picture. See \cite{GMM2024} for details.  

\section{Conclusion}\label{sect:conc}
In this paper we have shown that there are natural conditions on the orbits of break points of piecewise linear circle homeomorphisms which imply that such a map is conjugate to a rigid rational rotation. This condition
also has consequences for the existence of absolutely continuous probability invariant measures, jumps in the derivative
at break points, smoothness of the conjugacy and boundedness of derivatives of iterates. All of these related properties have analogues for PWL circle homeomorphisms which are conjugate to irrational rotations, see \cite{AA2016,Adouani2012,Akhad2006,Dzh2006,H1979,Khanin2003,Khmelev2005,Liousse2005} for example. 

\begin{figure}[b!]
\centering
\includegraphics[width=8cm]{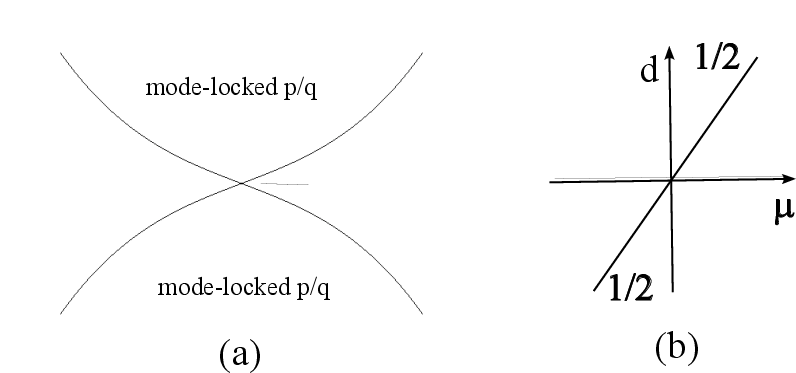}
\caption{
(a) Schematic view of a mode locked region in two parameter space showing pinching at the point at which the map is 
conjugate to a rigid rational rotation. (b) Linear approximation of the equivalent regions with $\frac{p}{q}=\frac{1}{2}$ for the example (\ref{ex:Hmod}) with $0<\lambda <1$. The rotation number is $\frac{1}{2}$ in the wedge-shaped regions labelled $\frac{1}{2}$.
}
\label{fig:pinch}
\end{figure} 

We have also considered this phenomenon in the context of families of maps. In stark contrast to the smooth case, 
the rotation number as a function of the parameter is differentiable at a value for which the map is conjugate to a rigid 
rational rotation. Thus the rotation number scales linearly near such points. As pointed out in section~\ref{sect:family} 
this is in general a high codimension phenomenon, but if there are only two break points then it is codimension two,
and if there are symmetries or other constraints then it can become codimension one, as in the examples of sections~\ref{sect:Hex} and \ref{sect:GMMex}. 

Under appropriate monotonicity assumptions, there is no mode-locking at the special points in parameter space 
at which the map is conjugate to a rigid rational rotation. Hence in a typical two-parameter family of maps the mode-locked
region associated must be `pinched' at such a point as shown in Figure~\ref{fig:pinch}. Pinching (shrinking points in their nomenclature) of mode-locked 
regions has been found in two-dimensional PWL systems (the border collision normal form) \cite{S2023,SM2009}, and it 
would be interesting to know whether this is a general mechanism for pinching or whether there are different
types of pinching phenomena for mode-locking in PWL maps. A simple illustration of this can be obtained from the Herman model (\ref{def:R}). Include a small offset $d$ to the position of the interior break point, so the break point becomes
\begin{equation}\label{eq:defcmod}
c=\frac{1}{1+\lambda}+d, \quad |d|\ {\rm small},
\end{equation}
then the modified equations for the lift are       
\begin{equation}\label{ex:Hmod}
x_{n+1} =\begin{cases}\mu +\frac{1}{1+\lambda}+\lambda x_n,  & {\rm if}~0\le x\le c\\ a+bx_n & {\rm if}~c<x<1\end{cases},
\end{equation}
where $\lambda$ is fixed and $a$ and $b$ are chosen so that the map is continuous at the break points. The equations for $a$ and $b$ are unhelpful, but if we work to order $\mu$ and $d$ (both assumed small) then correct to first order
\begin{equation}\label{eq:ab}
a\sim \mu+\frac{1}{1+\lambda}+1-\frac{1}{\lambda}-\frac{(1+\lambda)(1-\lambda^2)}{\lambda^2}d, \quad b\sim \frac{1}{\lambda}+\frac{(1+\lambda)(1-\lambda^2)}{\lambda^2}d.
\end{equation}
The boundary bifurcations creating orbits of period two occur when the break points are periodic, and these are straightforward to calculate to lowest order. The break point at $x=0$ maps to $\mu +\frac{1}{1+\lambda}$ which is greater than $c$ if $\mu>d$ and less than $c$ if $\mu<d$. If $\mu<d$ then we apply the map in $x<c$ again and the criterion for a periodic orbit on the circle is thus
\[
\mu +\frac{1}{1+\lambda}+\lambda\left( \mu +\frac{1}{1+\lambda} \right)=1, 
 \]
with solution $\mu =0$ ($d>0$). A similar but more complicated set of manipulations shows that the four boundaries are (to first order), with $F$ given by the right hand side of (\ref{ex:Hmod}) and its extension as a lift, $F(x+1)=F(x)+1$,
\begin{equation}
\begin{array}{ll}
F(0)<c,\ F^2(0)=1,   & \mu=0, \ d>0,\\
F(0)>c,\ F^2(0)=1,   & \mu=(1-\lambda )d, \ d<0,\\
F(c)<1,\ F^2(c)=c+1,  & \mu=0, \ d<0,\\
F(c)>1,\ F^2(c)=c+1, \ \ \ \ \ \ \ \  & \mu=(1-\lambda )d, \ d>0.
\end{array}
\label{eq:cases4}
\end{equation} 
Note that the curves are tangential at $(0,0)$ even though the maps are piecewise linear and that the curves $\mu =(1-\lambda )d$ change sign depending on the sign of $1-\lambda$ ($\lambda =1$ corresponds to a rigid rotation). Thus for example of $\lambda <1$ then the period two orbit exists (with rotation number $\frac{1}{2}$) in mode locked regions $0\le \mu\le (1-\lambda )d$ if $d>0$ and $(1-\lambda )d\le \mu\le 0$ if $d<0$. This is what is shown in Figure~\ref{fig:pinch}b.

\end{document}